\documentclass{article}
\usepackage{amsmath, amssymb, mathrsfs, eucal, latexsym}
\usepackage{color}
\usepackage{soul}
\usepackage{authblk}
\usepackage{graphicx}
\usepackage{geometry}

\newtheorem{theorem}{Theorem}[section]
\newtheorem{prop}[theorem]{Proposition}

\newtheorem{definition}[theorem]{Definition}
\newtheorem{example}[theorem]{Example}

\def\eps{\varepsilon}

\begin{document}

\title{Asymptotic hypothesis testing for the colour blind problem \\
}
\date{\today}
\author{Laura Dumitrescu } 
\author{Estate V. Khmaladze}
\affil{Victoria University of Wellington, New Zealand}

\maketitle
\begin{abstract}
\noindent In the classical two-sample problem, the conventional approach for testing distributions equality is based on the difference between the two marginal empirical distribution functions, whereas a test for independence is based on the contrast between the bivariate and the product of the marginal empirical distribution functions. In this article we consider the problem of testing independence and distributions equality when the observer is ``colour blind'' so he cannot distinguish the distribution which has generated each of the two measurements. Within a nonparametric framework, we propose an empirical process for this problem and find the linear statistic which is asymptotically optimal for testing the equality of the marginal distributions against a specific form of contiguous alternatives.

\vspace{1mm}

\end{abstract}

\noindent {\em AMS classification}: Primary 62G10; Secondary 62G20

\noindent {\em Keywords:} Asymptotically optimal test, contiguous alternatives, empirical process, goodness-of-fit, Kolmogorov-Smirnov statistics.

\section{Introduction} 
\label{introduction}
Consider an experiment when one observes pairs of balls with random diameters $\{(X_i, Y_i)\}_{1 \le i \le n}.$ The pairs are independent, identically distributed (i.i.d.) and the random diameters within each pair are also independent. The balls are coloured and the ball with diameter $X_i$ is green, while the other ball with diameter $Y_i$ is blue. We want to test if the diameters of the green balls and those of the blue balls have the same distribution. Denoting the cumulative   distribution functions by $P_1$ and $P_2,$ we want therefore to test the non-parametric hypothesis 
$$H_0: \ P_1=P_2 \;(\text {and equal to some unspecified} \; Q).$$
As it is formulated so far, the problem is a classical one and the class of test statistics with distributions independent of the common distribution $Q$ (provided it is continuous) is well known. Namely, if $P_{1n}$ and $P_{2n}$ are the empirical distribution functions (e.d.f.-s) of $\{X_i\}_{1 \le i \le n}$ and $\{Y_i\}_{1 \le i \le n},$ respectively, any statistic based on the empirical processes
$$ \sqrt{n}(P_{1n} - P_{2n}), \mbox { or } \sqrt{n} \left[P_{1n} - \frac{1}{2}(P_{1n} + P_{2n}) \right],$$
which are invariant under Kolmogorov time transformation $Q(x)=t,$ has a distribution which is the same for all $Q$. If we are interested only in the change between the expected values of the  diameters of the green and blue balls, then the Student's statistic $\sqrt{n}(\bar{X}_n - \bar{Y}_n)$ (with proper normalization) will provide a good test. 

However, what can one do in the case when an observer cannot distinguish the colours of the balls, that is, when he is ``colour blind''? 

In this case, it is not possible to construct the e.d.f.-s $P_{1n}$ and $P_{2n},$ or, even the averages $\bar{X}_n$ and $\bar{Y}_n.$ Yet, this paper shows that a systematic approach to testing $H_0$ is possible, again providing an empirical process and a test statistics, based on this process, with distribution independent of the (unknown) common $Q,$ and explains what is the price one pays for  ``colour blindness'' in terms of the power of our tests.
\vspace{2mm}

An example of the colour blind situation is encountered in the case of double-blind
trials, where comparing the effects of two treatments (e.g. placebo vs. a drug,
or an established drug vs. a new drug), is of interest. If each participant receives
each treatment once, at sufficiently distant moments in time, we may assume that
the independence between effects holds. By comparing the two treatment effects on each participant, the subject specific effect may be removed. Moreover, to eliminate subjectiveness in the evaluation of treatments, neither the observer/experimenter nor the participant can distinguish the treatment which is being administered, and only their effects are recorded. Double-blind studies are said to give more accurate results due to the potential reduction in the observer's/experimenter's bias.

Another example comes from genetics. In each nucleus of a somatic human cell, there are 23 pairs of chromosomes. Within each such pair, one chromosome is derived from the mother DNA and the other is derived from the father DNA. In karyotype analysis, measurements of different characteristics (such as the spiralization coefficient) are collected on homologous chromosomes, and the question of interest is to determine if there exist significant differences between the chromosomes derived from the mother and the chromosomes derived from the father. However, visually the chromosomes in the pair are not distinguishable. This example was the motivation behind the research presented in \cite{ParKhm82}. A detailed description of genetic data appears, in e.g. \cite{Tho00}.

%Finally, in \cite{cw-gol-kowa-ple82}, another example of observed pairs, whose order is unobservable was provided. Pairs of generation times of sister cells that result from the division of {\em Chilodonella steini} were analysed. In this case, a mother cell divides into two new cells: the anterior ({\em proter}) and the posterior ({\em opisthe}), but it may be difficult to distinguish between the two.

When observations are assumed to be collected from two independent normal populations, the estimation of means of unordered pairs of observations was considered in \cite{Hin73},  and for several populations, in \cite{BerSid72}.

A general case was proposed in \cite{ParKhm82}, where it was assumed that the two distributions $P_1$ and $P_2$  belong to the same parametric family of distributions, and a statistic for a locally most powerful rank test was derived. It was also shown that the degree of separation which can make the alternative hypothesis distinguishable from  the null, is of order $n^{-1/4}.$ As it will be seen in Section \ref{localalt}, an analogous finding is encountered in this article.

Under the assumption that the two observations within pairs are independent and belong to the family of Lehmann alternatives, a test for verifying their equality has been discussed in \cite{DaPh88}. A modified Mann-Whitney test statistic was constructed by taking into account the number of times the minimum in a pair exceeded the maximum from another pair, see also Section \ref{justify}. More recently, under the exponential tilting model for the two densities, and assuming, again, independence, a test based on the empirical Shannon information has been proposed in \cite{LiQ11}. 
\vskip 0.3cm

This article is organised as follows. In Section \ref{justify} we give a heuristic justification of the form of the empirical process and the test statistics which are introduced in Section \ref{ep}. In Section \ref{test_statistics} we show the behaviour of the Kolmogorov--Smirnov goodness of fit statistic and derive linear statistics, which are optimal for a particular sequence of local alternatives. Finally, in Section \ref{max}, we consider an empirical process based on the largest observation within each pair. This could be thought of as the first and the most direct object to consider. However, we think that  the roundabout way through the empirical process with two-dimensional time of the previous sections is actually simpler and more natural.

\section{Possibilities for testing}
\label{justify}
The only data that a colour blind observer can collect is the sequence of pairs $\{(U_{i}, V_{i})\}_{1 \le i \le n},$ where $U_{i}=\max\{X_i,Y_i\},$ and $V_{i}=\min\{X_i,Y_i\}.$ So, in  search for a test statistic, we investigate the relation between the distributions of the maximum and the minimum of each observed pair of independent observations. Note that regardless of whether the two marginal distributions are equal or not, the random variables $\{U_{i}\}_{1 \le i \le n}$ and $\{V_{i}\}_{1 \le i \le n}$ form two sequences of i.i.d.  random variables. Their cumulative distribution functions are, respectively, for $x \in \mathbb{R}$
\begin{eqnarray*}
	P^{(1)}(x) &=& P_1(x) + P_2(x) - P_1(x)P_2(x),\\
	P^{(2)}(x) &=& P_1(x) P_2(x).
\end{eqnarray*}
As a first important point,  note that, under $H_0$, the distributions $P^{(2)}$ and $P^{(1)}$ cannot be arbitrarily different; they are tied by the relation $P^{(2)}(x) = \left[1-\sqrt{1-P^{(1)}(x)}\right]^2.$ For arbitrary distributions $P_1$ and $P_2,$ the following inequalities hold, for any $x \in \mathbb{R},$ 
\begin{equation}
P^{(2)}(x) \le \left[ \frac{P^{(2)}(x)+P^{(1)}(x)}{2}\right]^2 \le \left[1-\sqrt{1-P^{(1)}(x)} \right]^2, \label{ineq-compare}
\end{equation} 
or, equivalently, $$P_1(x) P_2(x) \le \left[ \frac{P_{1}(x)+P_{2}(x)}{2}\right]^2 \le \left[1-\sqrt{(1-P_{1}(x))(1-P_{2}(x))} \right]^2,$$
with equality if and only if $H_0$ is true. %These relations show that, in the absence of a parametric form in $H_0,$  a proxy for the unknown common distribution can be taken to be $\displaystyle{Q=\frac{P_1+P_2}{2} = \frac{P_{(1)} + P_{(2)}}{2}.}$

\begin{figure}[h!]
\centering
\includegraphics[scale=0.6]{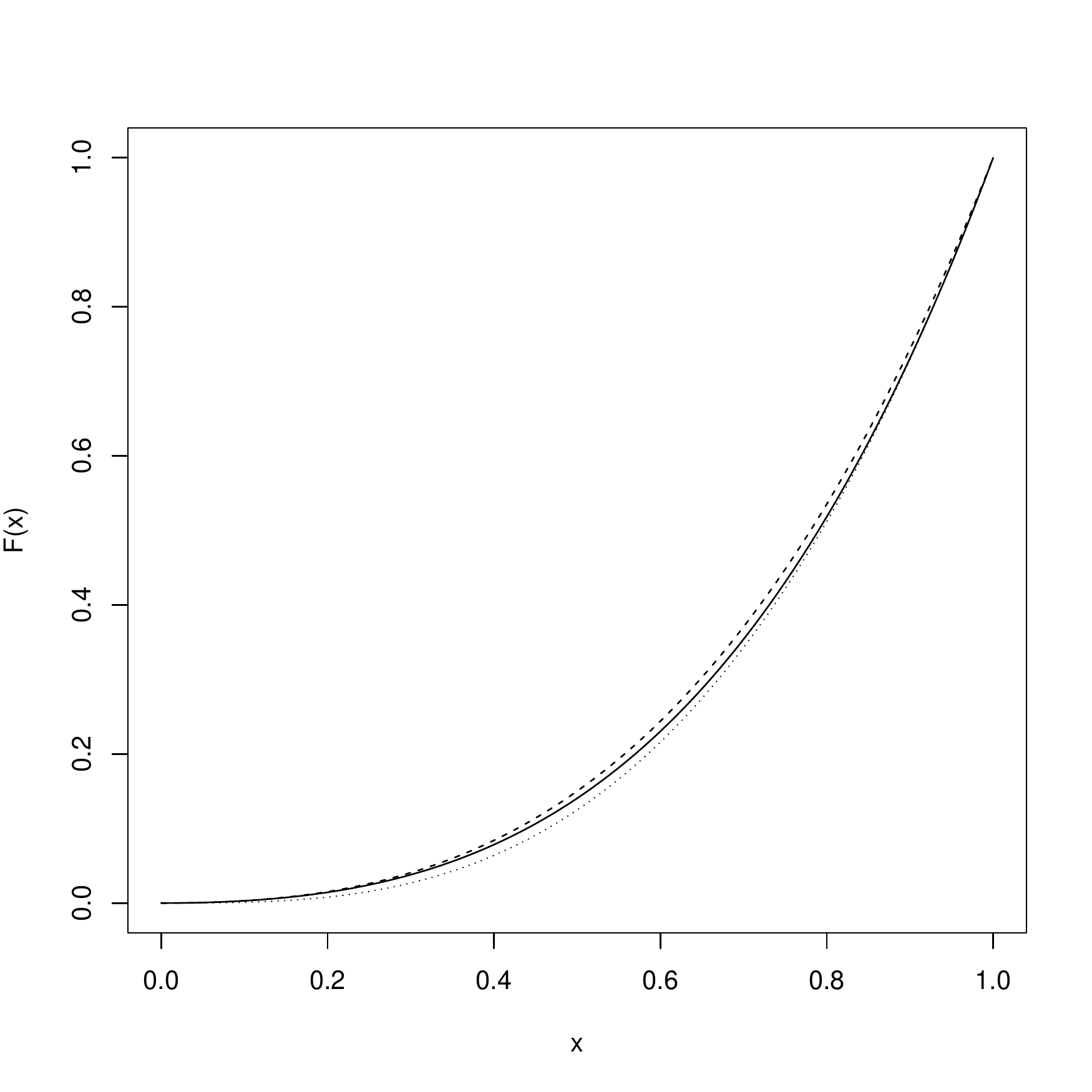}
\label{fig:graphs}
\caption{The case when $P_1(x)=x$ and $P_2(x)=x^2,$ with $0 \le x \le 1.$ The top dotted curve shows the graph of $[1-\sqrt{1-P^{(1)}(x)}]^2,$ the solid curve is a graph of $[(P_1(x)+P_2(x))/2]^2$, whereas the bottom dotted curve  gives the plot $P^{(2)}(x)$ . }
\end{figure}

As a next point, note that the inequalities appearing in \eqref{ineq-compare} are surprisingly tight even when the difference between $P_1$ and $P_2$ is considerable. Figure \ref{fig:graphs} illustrates this property for 
$P_1(x)=x$ and $P_2(x)=x^2$.
From a somewhat different side, observe the following fact although it is obvious that the maximum value within each pair will be greater than the minimum value in the same pair, comparing cross values will be informative. More precisely, if $i \neq j,$ under $H_0,$
\begin{eqnarray*}
P(U_j>V_i) = P(Q(U_j)>Q(V_i)) = \frac{5}{6},	
\end{eqnarray*}
while this probability would be larger if, for example, $P_1$ is stochastically dominated by $P_2,$ see \cite{ParKhm82} and \cite{DaPh88}.

These arguments suggest that a non-parametric testing approach in the colour-blind problem is possible, but deviations from $H_0$ will be difficult to detect. They also suggest that using a statistic which is based on an an empirical version of a contrast of the form  
$$\displaystyle{P_1 P_2 - \left(\frac{P_1+P_2}{2}\right)^2}$$
will be a natural approach. 

However, we prefer to %follow what may look somewhat roundabout route and 
first study the joint behaviour of $U_i$ and $V_i$. This allows, in our view, a natural treatment of the problem and helps to clarify the situation in Section \ref{max} below.\\

\section{The two-sample problem. The colour-blind process}
\label{ep}
In the sequel, without loss of generality, we assume that all observations $(X_1, Y_1), \ldots, (X_n, Y_n)$ are located in $[0,1]^2.$ We consider a more general approach to testing $H_0$ by using the empirical distribution of the pairs $\{(X_i, Y_i)\}_{1 \le i \le n},$ which are assumed to be i.i.d. two-dimensional random vectors, with a common bivariate distribution $\mathbb{Q}.$ More, precisely, for each $n \ge 1$ and $0 \le x \le 1, \ 0 \le y \le 1,$ let
$$\mathbb{P}_n(x,y)= \frac{1}{n}\sum_{i=1}^n {\bf 1}_{\{(X_i, Y_i) \in [0,x] \times [0,y]\}},$$ 
and define the product measure $\mathbb{Q}_n = Q_n \times Q_n,$ with marginals given by $$\displaystyle{Q_n = \frac{P_{1n} + P_{2n}}{2}} .$$ Here $P_{1n}(x) = 1/n\sum_{i=1}^n {\bf 1}_{\{X_i \le x\}}$ and $P_{2n}(x)=1/n \sum_{i=1}^n {\bf 1}_{\{Y_i \le x\}}$ are the empirical marginal distributions of $X$ and $Y,$ respectively but note that they cannot be obtained in the colour-blind problem. However, the distribution $Q_n$ can be computed as $P_{1n} + P_{2n}= P^{(2)}_n + P^{(1)}_n,$ where $P^{(2)}_n$ and $P^{(1)}_n$ are the corresponding empirical distribution functions based on the maxima $\{U_{i}\}_{1 \le i \le n}$ and the minima $\{V_{i}\}_{1 \le i \le n},$ respectively. 

Prompted by the discussion in Section \ref{introduction}, we consider the empirical process given by 
\begin{equation}\label{EP-not-CB}
\mathbb{R}_n(x,y)=\sqrt{n} \left[\mathbb{P}_n(x,y) - \mathbb{Q}_n(x,y) \right], \ 0 \le x \le 1, \ 0 \le y \le 1. 
\end{equation} 
Again, the empirical process $\mathbb{R}_n(x,y)$ can not be used  directly because the empirical distribution function $\mathbb{P}_n,$ can not be obtained in the colour-blind problem. However, this distribution can be obtained within the class of symmetric Borel sets, i.e., sets $B$ for which  $(x,y)\in B$ implies $(y,x)\in B$. Indeed, if $B$ is symmetric, then for any $i \le n,$ we have
$$P((U_{i}, V_{i}) \in B) = P((X_{i}, Y_{i}) \in B)$$
and a similar equality holds for the empirical distribution functions. Therefore, we consider the process $\mathbb{R}_n$ on symmetric sets, as in the following definition. Denote by $S_{x,y}=[0,x]\times [0,y]\cup [0,y]\times [0,x]$, which gives the symmetrised version of the rectangle $[0,x]\times [0,y].$ By construction, if $u=\max\{x,y\}$ and $v=\min
\{x,y\}$ then $S_{u,v}=S_{v,u}$.

\begin{definition}\label{def2.1}  Let $\cal{B}$ denote the class of symmetric Borel subsets of $[0,1]^2$. Then, the restriction of $\mathbb{R}_n$ to  $\cal{B}$ is called a colour-blind empirical process
$$ \mathbb{R}_n^{s}(B) = \mathbb{R}_n(B), \; \text{where} \; B\in \cal B.$$
For $S_{u,v}$ as defined above, we use the notation 
\begin{align*}
\mathbb{R}_n^{s}(u,v) = \mathbb{R}_n(S_{u,v}) = \mathbb{R}_n(u, v)+\mathbb{R}_n(v, u)-\mathbb{R}_n(v,v),\end{align*}
\end{definition}	
and so, essentially, the process $\mathbb{R}_n^{s} (u,v)$ is defined on the simplex $0\le v\le u\le1$.

The functional central limit theorem for the empirical process $\mathbb{R}_n$ holds without the symmetry assumptions. Indeed, if $v_n$ denote the  ``usual" empirical process, based on i.i.d. pairs $\{(X_i, Y_i)\}_{1 \le i \le n},$ with common distribution $\mathbb{Q},$
$$v_n(x,y) = \sqrt{n}(\mathbb{P}_n - \mathbb{Q})(x, y) ,$$ 
and if $v_{\mathbb{Q}}$ denotes the Brownian bridge in time $\mathbb{Q}$, then, cf. \cite{vdVW96},
$$v_n \stackrel{\mathcal{D}} \longrightarrow v_{\mathbb{Q}}, \; \mbox{as } n\to \infty.$$
However, under $H_0,$ $\mathbb{Q}(x,y)= Q(x) Q(y)$ and so, the empirical process appearing in \eqref{EP-not-CB} can be written as
\begin{align} \label{process-NCB}
\mathbb{R}_n(x,y) &= v_n(x,y) - \sqrt{n}[Q_n(x)Q_n(y) - Q(x)Q(y)] \nonumber \\
 	&= v_n(x,y) - \sqrt{n}[Q_n(x) - Q(x)]Q(y) - \sqrt{n}[Q_n(y) - Q(y)]Q(x) + r_n(x,y) \nonumber\\
 	&= v_n(x,y) - \frac{1}{2}[v_n(x,1)+v_n(1,x)] Q(y) - \frac{1}{2}[v_n(1,y)+v_n(y,1)] Q(x) + r_n(x,y), 
 \end{align}
 where 
 \begin{equation}\label{small}
 r_n(x,y) = \sqrt{n} [Q_n(x) -Q(x)][Q_n(y) - Q(y)],
 \end{equation}
and $\sup_{(x,y)\in [0,1]^2} |r_n(x,y)| = o_P(1).$

Since the leading term in the right hand side of \eqref{process-NCB} is a linear transformation of $v_n$, the central limit theorem for $\mathbb{R}_n$ and $\mathbb{R}_n^s$ follows and so, as $n\to\infty$, we have $\mathbb{R}_n \stackrel{\mathcal{D}} \longrightarrow \mathbb{R}$, where
\begin{equation}\label{proj}
\mathbb{R}(x,y) = v_{\mathbb{Q}}(x,y) - \frac{1}{2}[v_{\mathbb{Q}}(x,1)+v_{\mathbb{Q}}(1,x)] Q(y) - \frac{1}{2}[v_{\mathbb{Q}}(1,y)+v_{\mathbb{Q}}(y,1)] Q(x).
\end{equation}
It follows that $\mathbb{R}_n^s\stackrel{\mathcal{D}} \longrightarrow \mathbb{R}^s$, where $\mathbb{R}^s$ is the restriction of the process $\mathbb{R}$ to the set $\cal B$.

The transformation of $v_n$ in the right hand side of \eqref{process-NCB}  (as well as the transformation of $v_{\mathbb{Q}}$ in \eqref{proj}) is, actually, a projection. However, for better insight in the nature of the process $\mathbb{R}_n^s$ let us consider another projection of $v_{\mathbb{Q}}$ given by the  operator
$$\mathcal{L }\alpha(x,y) = \alpha(x,y) - Q(x)\alpha(1,y)  -Q(y)\alpha(x,1)  + Q(x)Q(y) \alpha(1,1) .$$
It projects a function $\alpha$ onto the class of functions equal to zero  everywhere on the boundary of $[0,1]^2$. By applying it to $v_{\mathbb{Q}}$ we obtain
\begin{equation}
z_{\mathbb{Q}}(x,y) = \mathcal{L }v_{\mathbb{Q}} (x,y) = v_{\mathbb{Q}}(x,y) - Q(x) v_{\mathbb{Q}}(1,y) - Q(y)v_{\mathbb{Q}}(x,1)  + Q(x) Q(y) v_{\mathbb{Q}}(1,1)\label{pillow-def}
\end{equation}
as a projection of $v_{\mathbb{Q}}$. The process $z_{\mathbb{Q}}$ is called a {\em Brownian pillow} on $[0, 1]^2$ (or a {\em bivariate tied-down Brownian bridge} or {\em completely tucked Brownian sheet}) and it appears as the limit process for  testing independence of components of continuous bivariate random vectors, based on empirical distributions functions (see e.g. \cite{BlKieRos61} and the Section 3.8 in \cite{vdVW96}).
 By time transformation $t=Q(x), \ s=Q(y),$ it can be mapped into a {\it standard} Brownian pillow in $t$ and $s$, i.e. a Gaussian process $z(s,t)$, with continuous sample paths and covariance function given by 
$ {\rm E} (z(s', t') z(s'', t'')) = (\min\{s', s''\} - s's'')(\min\{t', t''\} - t't''), \ 0 \le s',s'',t',t'' \le 1.$
 Its finite $n$ version is, obviously, given by $z_n=\mathcal{L} v_n$.

Our interest in $z_{\mathbb{Q}}$ and $z_n$ stems from the following fact.

\begin{prop}\label{sym-pillow} For symmetric sets $B\in\cal B$ we have 
$$ \mathbb{R}_n^s (B) = z_{n}(B) + o_P(1), \; \mbox{as } n\to\infty, \; \text{and} \quad  \mathbb{R}^s (B) = z_{\mathbb{Q}}(B) .$$
In particular, for any $n \ge 1$ and symmetrised rectangles, the following relationship holds 
$$\mathbb{R}_n^s  (u,v) = z_n(S_{u,v}) + r_n(S_{u,v}),$$ 
where 
\begin{eqnarray*}
z_n(S_{u,v}) &=&  z_n(u, v)+z_n(v, u)-z_n(v,v), \\
r_n(S_{u,v}) &=&  2\sqrt{n} [Q_n(u) -Q(u)][Q_n(v) - Q(v)] - \sqrt{n} [Q_n(v) -Q(v)]^2 =o_P(1).
\end{eqnarray*}
\end{prop}
Thus, the proposition describes the limit in distribution of the colour-blind process as a restriction of Brownian pillow to the class of symmetric sets. 

\subsection{Description of local alternatives}
\label{localalt}
We are interested in detecting small departures from the null hypothesis $H_0$ and, assuming independence between the two marginal distributions, such deviations will be specified by a sequence of probability distributions of the form $\{A_{1n} \times A_{2n}\}_{n \ge 1}$. 

Let $Q$ be an arbitrary, but fixed probability distribution. Let each $A_{1n}$ and $A_{2n}$ be two continuous probability distributions which are defined as  asymptotically ``small" departures from $\mathbb{Q}.$ 
%\begin{eqnarray*}
\begin{equation}\label{RN}
 \frac{dA_{1n}}{dQ} (x) = 1 + \eps_n h_{1n}(x), \quad
\frac{dA_{2n}}{dQ} (x) = 1 + \eps_n h_{2n}(x),
%\end{eqnarray*}
\end{equation}
where $\eps_n\to 0$ as $n\to\infty$, and the functions $h_{kn}, k=1,2$ converge to square integrable functions $h_{k} (x)$
$$\int_0^1 [h_{kn} (x) - h_k(x)]^2 dQ(x) \to 0, \ \int_0^1 h_k^2(x) dQ(x) < \infty, \ k=1,2.$$

From the definition of $h_{kn}$ it follows that, for all $n \ge 1$, $\int_0^1 h_{kn} (x) dQ(x) = 0, \ k=1,2$, and their limits inherit this property. 

\begin{figure}[h!]
\centering
\includegraphics[scale=0.6]{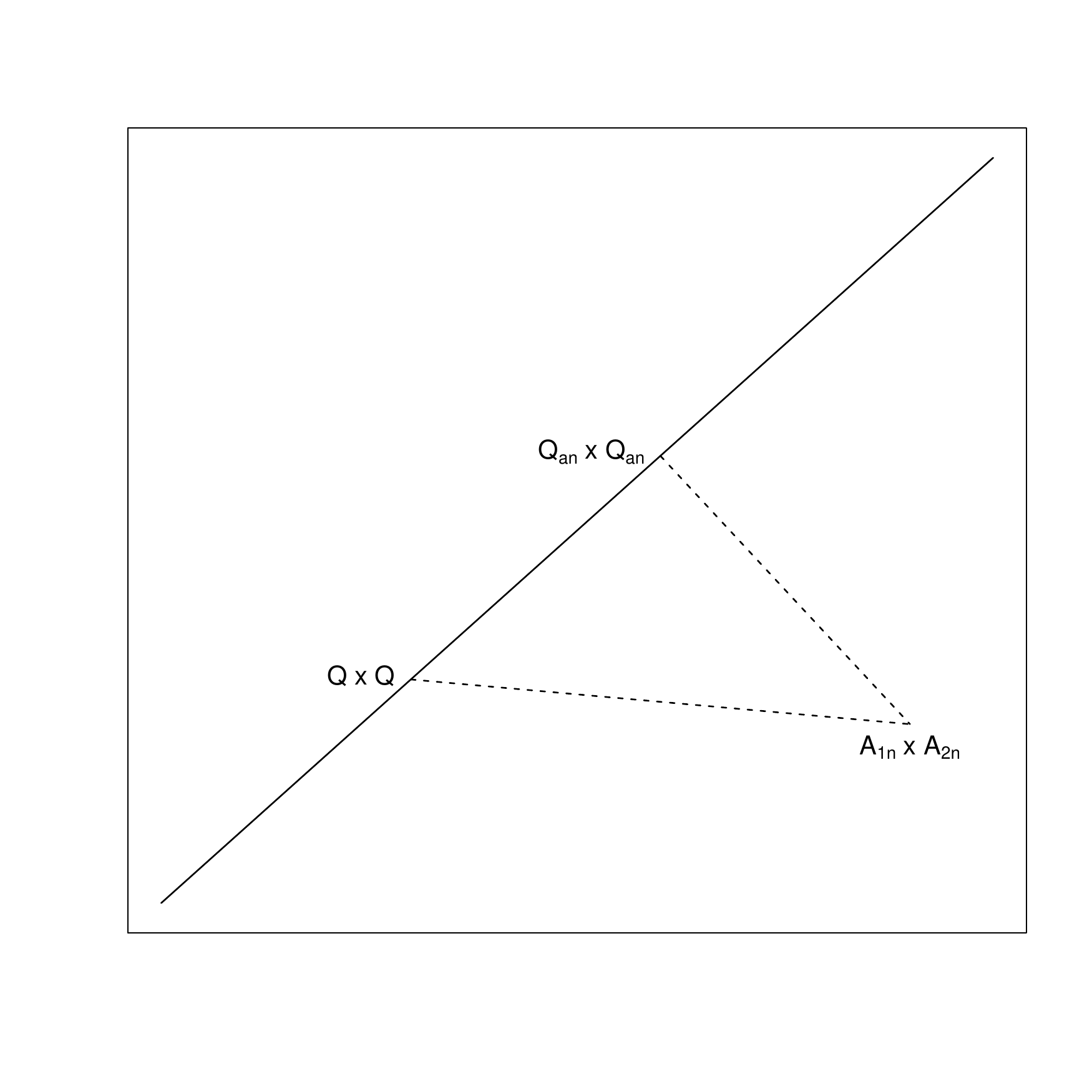}
\caption{The alternative $A_{1n}\times A_{2n}$ may be some distance away from the hypothetical pair $Q\times Q$, to which it converges, but statistics form the empirical process will ``react" on this alternative as much as it deviates from its projection $Q_{an}\times Q_{an}$.}
\label{fig:alts}
\end{figure}

Later we will see that $\eps_n$ will not be of order $1/\sqrt{n}$ as it is typically the case within the theory of contiguity \cite{OvZ12}, but will need to decrease slower. Therefore, we will eventually be outside the contiguity theory, and, therefore, we can neglect using square roots from the Radon - Nikodym  derivatives above, cf. \cite{OvZ12}, \cite{Rou72}, and consider them as they are, which is somewhat simpler. It is more interesting to recall that although both $A_{1n}$ and $A_{2n}$ tend to $Q$, in testing $H_0$ they will look differently. To see this fact, specific to the two-sample problem (and not to colour-blindness as such), consider the expected value of the classical two-sample process $\sqrt{n} (P_{1n} - P_{2n})$ under the alternative $\mathbb{A}_n=A_{1n} \times A_{2n}$. Introduce the functions $H_{1n}$ and $H_{2n}$ as
$$H_{kn}(x)=\int_0^x h_{kn}(y)dQ(y), \; \text{so that}
\; \; H_{kn}(0)=H_{kn}(1)=0, \; k=1,2,
$$
and so 
\begin{align*}
{\rm E}_a \sqrt{n} [P_{1n}(x) - P_{2n}(x)]&=\sqrt{n} \left[A_{1n} (x) - A_{2n} (x)\right ]\\
&= \sqrt{n}\eps_n [H_{1n}(x) - H_{2n}(x)]. 
\end{align*}
From this it can be shown that one can choose $\eps_n=1/\sqrt{n}$ and the linear statistic
$$\int_0^1 \phi(x) \sqrt{n} [dP_{1n}(x) - dP_{2n}(x)], \text{where} \; \; \phi(x) = h_{1n} (x) - h_{2n}(x) ,$$
leads to the asymptotically most powerful test, among those based on $\sqrt{n} (P_{1n} - P_{2n})$. However, the power of this test is only less than or equal to the power of the optimal (Neyman - Pearson) test in the problem of discriminating between alternative $\mathbb{A}_n$ and the hypothesis $Q\times Q$, cf. \cite{HaSiSe99}. To obtain an equality one has to change $Q$ to $Q_{an}=(A_{1n} + A_{2n})/2$ and note that this agrees with the form in Section \ref{justify}. Then, the test based on the linear statistic above becomes asymptotically equivalent to the Neyman - Pearson test for discriminating $\mathbb{A}_n$ and the hypothesis $Q_{an}\times Q_{an}$. The Radon-Nikodym derivatives of $A_{1n}$ and $A_{2n}$ with respect to $Q_{an}$ posses a symmetric form
\begin{equation*}
\frac{dA_{1n}}{dQ_{an}} (x) = 1 + \eps_n h_n(x), \quad
\frac{dA_{2n}}{dQ_{an}} (x) = 1 - \eps_n h_{n}(x),
%\end{eqnarray*}
\end{equation*}
with 
$$ h_n = \frac{(h_{1n} - h_{2n})/2}{1+\eps_n(h_{1n} + h_{2n})/2}.$$
Dependence of $Q_{an}$ on $n$, which itself converges to $Q$, is immaterial and we can assume from now on that in \eqref{RN} we have $h_{1n} = -h_{2n}$, i.e., we consider the class of local alternatives of the form
\begin{equation}\label{RN1}
\frac{dA_{1n}}{dQ} (x) = 1 + \eps_n h_{n}(x), \quad
\frac{dA_{2n}}{dQ} (x) = 1 - \eps_n h_{n}(x),
\end{equation}
for some continuous distribution $Q$ on $[0,1],$ assuming that there exists a function $h \in L^2(Q)$ such that
$$\int_0^1 h(x) dQ(x) =0, \text{and} \; \int_0^1[ h_n(x) - h(x)]^2 dQ(x) \to 0.$$ 
The Figure \ref{fig:alts} illustrates the situation.
The function $h$ determines the direction in which the alternative distribution $\mathbb{A}_n$ approaches the distribution $\mathbb{Q}.$

Now consider the rate of convergence of $\eps_n$ in the colour-blind problem. Under $H_a,$ 
\begin{eqnarray*}
{\rm E}_a[v_n(x,y)]&=&{\rm E}_a \{\sqrt{n}[\mathbb{Q}_{n}(x,y) - \mathbb{Q}(x,y)]\}\\
&=& \sqrt{n}[-\varepsilon_n Q(x) H_n(y) + \varepsilon_n H_n(x) Q(y) - \varepsilon_n^2 H_n(x) H_n(y)],
\end{eqnarray*}
where, as above, $H_n(x) = \int_0^x h_n(y) dQ(y).$ Since $H_n(1)=0,$ it follows that
\begin{eqnarray*}
{\rm E}_{a}[\mathbb{R}_n(x,y)] &=&  \sqrt{n}[-\varepsilon_n Q(x) H_n(y) + \varepsilon_n H_n(x) Q(y) - \varepsilon_n^2 H_n(x) H_n(y)],
\end{eqnarray*}
which shows that the statistics based on the process $\mathbb{R}_n(x,y)$ can distinguish alternatives with $\eps_n=O(n^{-1/2}).$	However, in the case of the symmetrized  process $\mathbb{R}_n^{s}(u,v),$ under $H_a,$ the linear term in $\eps_n$ becomes zero and we have
\begin{equation}
{\rm E}_{a}[\mathbb{R}_n^{s}(u,v)] =  \sqrt{n} \varepsilon_n^2[-2 H_n(u) H_n(v) + H_n(v)^2]. \nonumber
\end{equation}
This shows the loss of power when using the colour-blind statistic: only alternatives with $\varepsilon_n=O(n^{-1/4})$ can be detected by the colour-blind process.

%We conclude with a remark of a technical nature:  there is a subtle difference between conditions \eqref{RNa} with $h_1$ and $h_2$ being the same function and \eqref{RN1}, as the comparison with \cite{ParKhm82} reveals. The asymptotic expression for the optimal linear statistic in the Proposition \ref{as-opt} in the two cases may be slightly different, and simpler in this paper. This may happen as the result of possible influence of $\varepsilon^2$ terms for $\varepsilon\sim n^{-1/4}$, particularly in parametric models.

\subsection{Local alternatives of dependence}
\label{testing-independence}
It took some time to realise that it is also possible to test independence of diameters of the coloured balls in the colour-blind situation. The tests can be based on the same colour-blind process and the problem is actually easier than the problem of testing equality of distributions: it is possible to detect local alternatives with rate $1/\sqrt{n}.$

Introduce the {\em new class} of alternative hypotheses of the form
$$H_a': \ \mathbb{A}_n' = \mathbb{Q} + \varepsilon_n \mathbb{G}_n,$$
\begin{align}
\text{where } &\frac{d\mathbb{A}_n'(x,y)}{d\mathbb{Q}(x,y)}=1+\varepsilon_n g_n(x,y), \ g_n\in  L^2(Q\times Q), \text{ and there exists } g \in  L^2(Q\times Q) \text { such that } \notag \\
 &\  \int_0^1 \int_0^1 [g_n(x,y) - g(x,y)]^2 dQ(x)dQ(y) \to 0, \ \mbox{as } n \to \infty. \notag
\end{align}
Again, the function $g(x,y)$ here describes the functional direction, from which the sequence of alternatives $\{\mathbb{A}'_n\}_{n \ge 1}$ approaches $\mathbb{Q}.$ In addition, we specify the following integral conditions on $g_n,$ assuring that the marginal distributions of $X_i$ and $Y_i$ under $\mathbb{A}_n'$ are the same as under $\mathbb{Q}$
\begin{align} \label{RN2}
\int_0^1 g_n(x,y) dQ(x) = \int_0^1 g_n(x,y) dQ(y) = 0.
\end{align}
In comparison, the local alternatives introduced in \eqref{RN1}, which are used for testing the identical distributions assumption, approach $\mathbb{Q}$ from the functional direction $h(x)-h(y).$ This difference, as a function in  $L^2(Q\times Q)$, is orthogonal to the function $g.$ Thus, alternatives \eqref{RN1} and \eqref{RN2} approach the hypothetical family from orthogonal directions. As the following paragraph shows, the neighbourhood of the hypothetical family in these different directions is quite non-homogeneous.

Indeed, with $\mathbb{G}_n(x,y)=\int_0^x\int_0^y g_n(x',y') Q(dx') Q(dy')$, one can see that
\begin{equation*}
{\rm E}_{a'}[v_n(x,y)]={\rm E}_{a'}\{ \sqrt{n}[\mathbb{P}_{n} (x,y)- \mathbb{Q}(x,y)]\} = \sqrt{n}\varepsilon_n \mathbb{G}_n(x,y),
\end{equation*}
and therefore, noting that $G_n(1,y)=G_n(x,1) = 0$, we also have
\begin{align*}
{\rm E}_{a'}[\mathbb{R}_n(x,y)] &=  \sqrt{n} \varepsilon_n G_n(x,y), \\  {\rm E}_{a'}[\mathbb{R}_n^{s}(u,v)] &= \sqrt{n}\varepsilon_n [G_n(u,v)+G_n(v,u)- G_n(v,v)].
\end{align*}
Therefore the ``usual" rate of $n^{-1/2}$ will render the shift ${\rm E}_{a'}\mathbb{R}_n^{s}$ positive.

However, the symmetrisation again can make an alternative undetectable. To see that, let  $$ g_n^*(x,y) = [g_n(x,y) - g_n(y,x)]/2$$
be the anti-symmetric part of $g_n.$ We remark that it will make zero contribution to the shift of ${\rm E}_{a'}\mathbb{R}_n^{s}$, so that if $g_n$ is itself anti-symmetric, then the shift of ${\rm E}_{a'}\mathbb{R}_n^{s}$ becomes $0$, and the alternatives $\mathbb{A}_n'$ become undetectable. However, we will not pursue this approach here any further.

\section{Test statistics}
\label{test_statistics}

\begin{figure}[h!]
\centering
\includegraphics[scale=0.6]{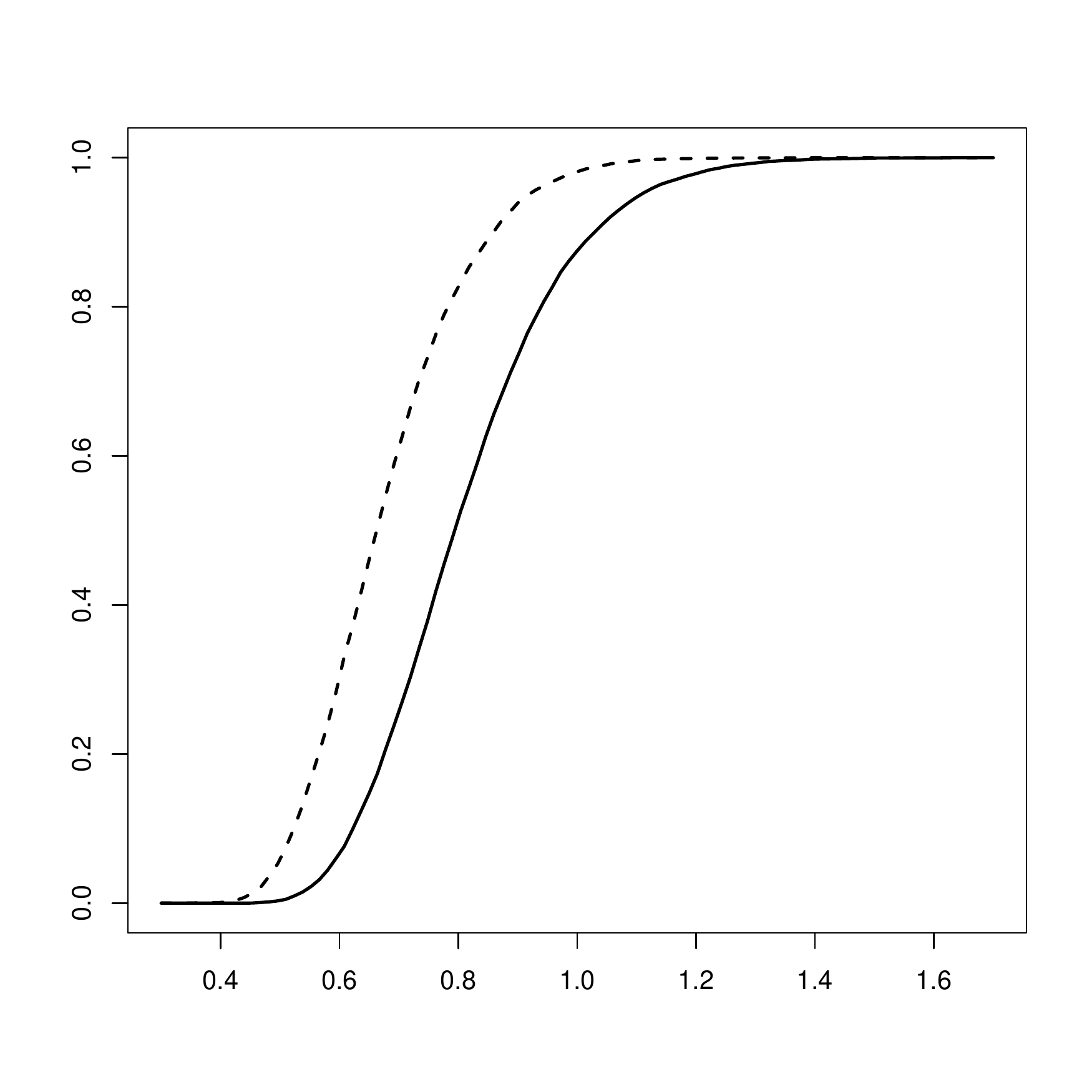}
\caption{The dotted line represents the graph of the (simulated) distribution function of $D_n^s,$ while the solid line represents that of $D_n$. The picture is consistent with the fact that $\mathbb{R}^s_n(u,v)$ is restriction of $\mathbb{R}_n(x,y)$, cf. Definition \ref{def2.1}. The size of the generated sample was $n=1,000$.}
\label{fig:SymandNot}
\end{figure}

\begin{figure}[h!]
\centering
\includegraphics[scale=0.6]{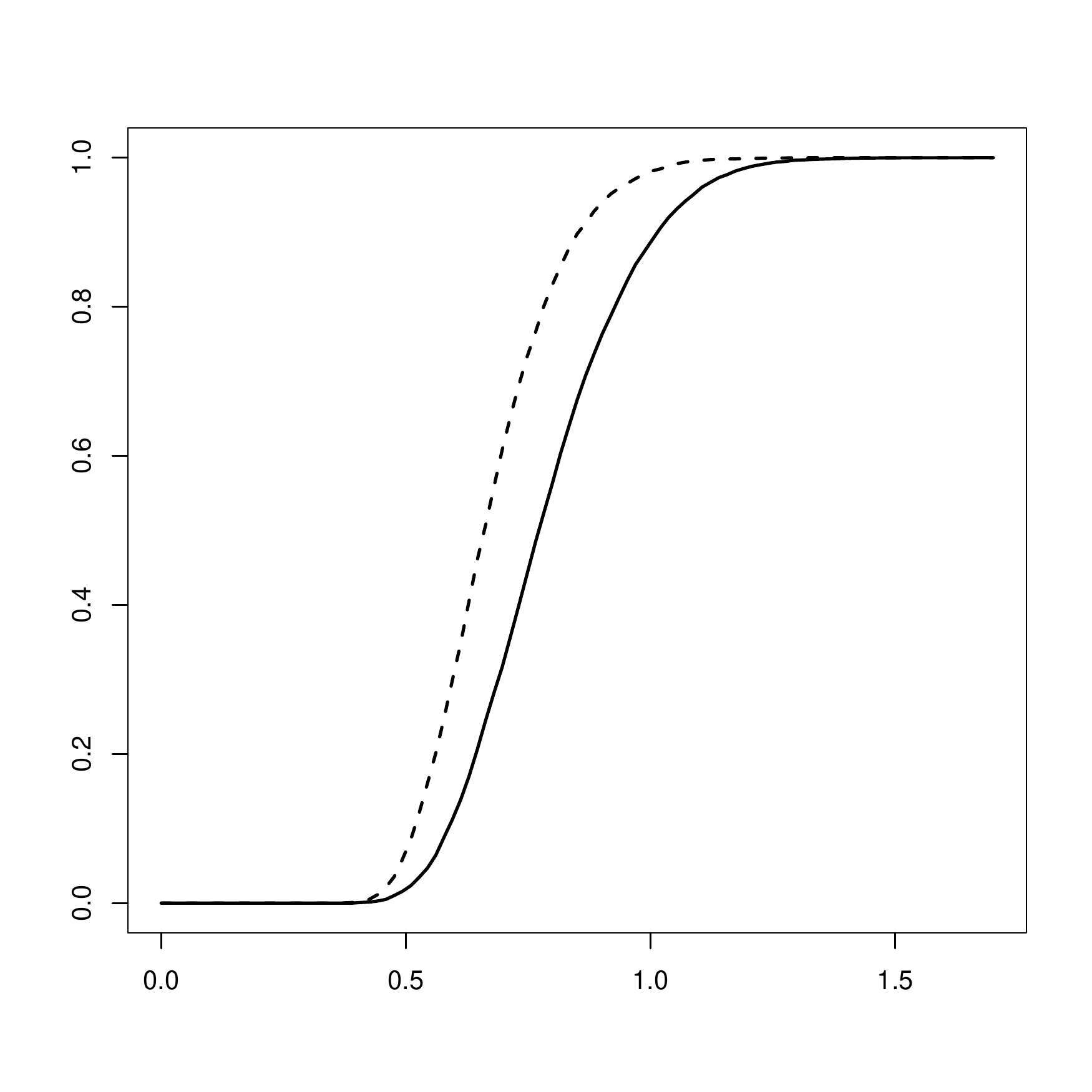}
\caption{The dotted line shows the simulated distribution function of $D_n^s$ under the null hypothesis, while the solid line shows its distribution under the alternative $A_1(x)=x, A_2(x)=x^2$. Figure \ref{fig:graphs} has illustrated that in the colour-blind situation, the alternative distributions look surprisingly difficult to distinguish from the null. However, according to the present graph, with a number of observations equal to $n=500$, the Kolmogorov-Smirnov test will have some power.}
\label{fig:KSalt1}
\end{figure}

\begin{figure}[h!]
\centering
\includegraphics[scale=0.7]{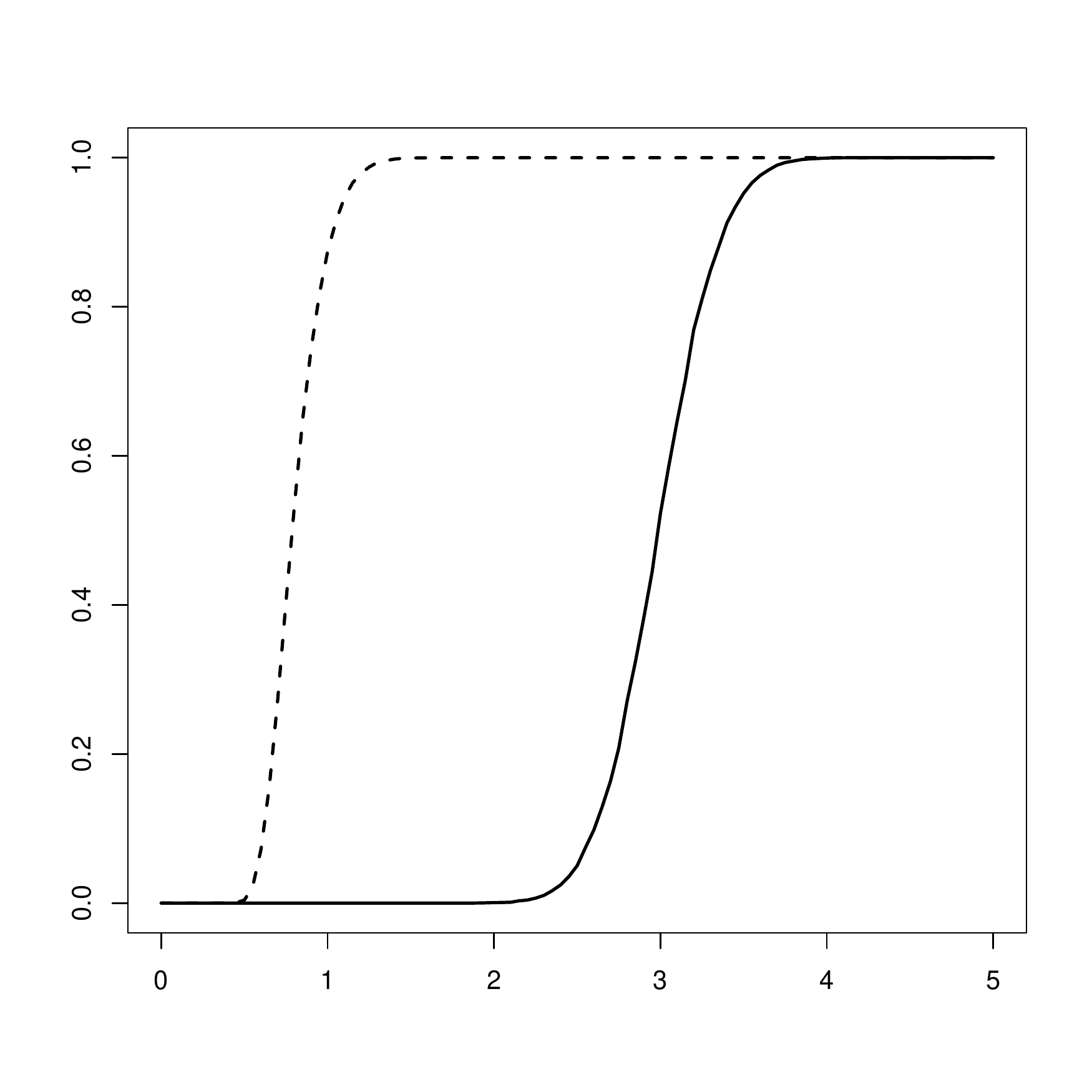}
\caption{The graph shows the simulated distribution function of $D_n$ under the null hypothesis (dotted line) and under the alternative $A_1(x)=x, A_2(x)=x^2$ (solid line). With a number of observations equal to $n=500,$ the discrimination between them would be obvious.}
\label{fig:KSalt2}
\end{figure}

Proposition \ref{sym-pillow} shows that, asymptotically, the colour-blind process $\mathbb{R}^s_n$ is equivalent to a Brownian pillow on symmetric sets. We also noted that the Brownian pillow $z_{\mathbb{Q}}$ can be transformed into standard Brownian pillow, i.e., a Brownian pillow when $\mathbb{Q}$ is the Lebesgue measure on $[0,1]^2$. Therefore, the classical goodness of fit statistics, such as the Kolmogorov-Smirnov statistic, and other statistics based on $\mathbb{R}^s_n$ and invariant with respect to time transformation $t=Q(x), s=Q(y)$ will a have limiting distribution, which does not depend on unknown $Q$. Intuitively, one would expect that among the two statistics 
$$D_n=\sup_{(x,y)\in [0, 1]^2} |\mathbb{R}_n(x,y)| ,\quad \; D_n^s = \sup_{(u,v)\in [0, 1]^2, v<u} |\mathbb{R}^s_n(u,v)| , $$ 
the second one will be stochastically smaller.  Figure \ref{fig:SymandNot} shows the graphs of their distributions and confirms that this intuition is correct. 
To illustrate the situation in terms of the power of goodness of fit tests, we consider in Figure \ref{fig:KSalt1} the shift of $D_n^s$ under the alternatives described earlier in Figure \ref{fig:graphs}.
Furthermore, to better illustrate the consequence of colour-blindness, in Figure \ref{fig:KSalt2}, we present the simulated distribution functions of $D_n,$ under the null and under the same alternatives as above, or as in Figure \ref{fig:graphs}. We observe that the discrimination between the two could have been absolutely obvious.

\subsection{Linear statistics}

Unlike the goodness of fit tests, which are of omnibus nature and typically have some power against a very wide class ot alternatives, the tests based on linear statistics may have asymptotically no power against the ``majority" of alternatives, but are asymptotically most powerful against a certain form of  alternatives, with a specific function $h_n$. In this section we derive the form of such statistic in the colour-blind problem.

Let $\varphi \in L^2(Q \times Q)$ and consider the function-parametric version of the colour-blind process introduced in Definition \ref{def2.1}, i.e. $$\mathbb{R}_n^{s}(\varphi) = \int_0^1\int_{v<u} \varphi (u,v)d\mathbb{R}_n^{s} (u,v) .$$ In search for the optimal linear functional, the next result shows that we can restrict our attention to symmetric functionals from $ z_n.$
\begin{prop} For every $\varphi \in L^2(Q \times Q)$, as $n \to \infty$,
\begin{equation}
 \int_0^1 \int_{v<u} \varphi(u,v) d\mathbb{R}_n^{s}(u,v) = \int_{0}^1 \int_{0}^1 \tilde{\varphi}(x, y) dz_n(x,y) + o_P(1), \label{fn-param-equiv}
\end{equation}
where \[ \tilde{\varphi}(x, y) = \begin{cases} 
\varphi(x,y), & x \geq  y \\
\varphi(y,x), & x < y 
\end{cases}.
\]
\end{prop}
{\bf Proof.} First note that the planar integral from the term $z_n(v,v)$ in Proposition \ref{sym-pillow} is zero and that the integral from the residual term is indeed small, i.e. for any $\varphi \in L^2(Q \times Q)$
$$\int_0^1\int_0^1 \varphi(x,y)  \sqrt{n} [dQ_n(x) -dQ(x)][dQ_n(y) - dQ(y)]  = o_P(1).$$
Then, due to Proposition \ref{sym-pillow}, the left hand side of \eqref{fn-param-equiv} can be written as
\begin{eqnarray}
\int_0^1\int_{v<u} \varphi(u,v) d\mathbb{R}_n^{s}(u,v) = \int_0^1\int_{v<u} \varphi(u, v) [d{z}_n(u,v)  +  d{z}_n(v,u)] + o_P(1).  \label{LHS-equiv-fn-param}
\end{eqnarray}
The main term in the right hand side of \eqref{fn-param-equiv} is 
\begin{eqnarray*}
\int_{0}^1 \int_{0}^1 \tilde{\varphi}(x,y) dz_n(x,y) 
&=& \int_0^1\int_{x >y} \varphi(x,y) dz_n(x,y)+ \int_0^1 \int_{y >x}\varphi(y,x) dz_n(x,y)\\
&=& \int_0^1 \int_{u >v} \varphi(u,v) dz_n(u,v)+\int_0^1 \int_{u >v}\varphi(u,v) dz_n(v,u),
\end{eqnarray*}
where we have used the symmetry property of $\tilde \varphi$. The proof is concluded by noting that the right hand sides of the last two displays are equal. \hfill $\Box$\\

In what follows, we assume that $\varphi \in L^2(Q \times Q)$ is a symmetric function. 

Denote by $\displaystyle{(Q\varphi)(x) = \int_{0}^1 \varphi(x,y) dQ(y)}$ and consider the projection of $\varphi$ given by 
$$(\mathcal{L}^* \varphi)(x,y)=  \varphi(x,y) - (Q\varphi)(x) - (Q\varphi)(y) + {\rm E}_{\mathbb{Q}} [\varphi(x,y)].$$ 

The process $z_n$ was introduced as a projection of $v_n$ and so $z_n = \mathcal{L} v_n,$ with $(\mathcal{L} v_n)(x,y) = v_n(x,y) -Q(x) v_n(1,y) - Q(y)v_n(x,1) + Q(x) Q(y) v_n(1,1).$ The following proposition shows that operator $\mathcal{L }^*$ can be viewed as the adjoint of the operator $\mathcal{L }$ introduced earlier, cf. \cite{Kuo75}.

\begin{prop}\label{z_n-functionals} We have the following  
\begin{align*} z_n(\varphi) &= \mathcal{L}v_n(\varphi) =\int_{0}^1 \int_{0}^1 \varphi (x,y) dz_n(x,y)\\
&= \int_{0}^1 \int_{0}^1 (\mathcal{L}^*\varphi)(x,y) dv_n(x,y)= v_n(\mathcal{L}^*\varphi) .
\end{align*}
\end{prop}
{\bf Proof.} By direct computation,
\begin{align*}
z_n(\varphi) &= \int_{0}^1 \int_{0}^1 \varphi (x,y) dv_n(x,y) - \int_{0}^1 \int_{0}^1 \varphi (x,y) dQ(x) dv_n(1,y) \\
&- \int_{0}^1 \int_{0}^1\varphi (x,y) dQ(y) dv_n(x,1)\\
&= \int_{0}^1 \int_{0}^1 \left[\varphi(x,y) - \int_{0}^1 \varphi(x,y) dQ(x) - \int_{0}^1 \varphi(x,y) dQ(y) \right] dv_n(x,y). \; \hfill \Box
\end{align*} 
As a result, \begin{equation}
\label{fn-parametric-equiv}
\mathbb{R}_n^{s}(\varphi) = v_n(\mathcal{L^*} \tilde{\varphi}) + o_P(1) ,
\end{equation}
and we can use the well-known central limit theorem for the function-parametric empirical process $v_n$ to describe the asymptotic behaviour of  $\mathbb{R}_n^s(\varphi)$.

\subsection{Optimal linear statistics}
\label{Opt_lin}
Under $H_a,$ the expected value of $v_n(\tilde{\varphi}) $ is not zero. For every symmetric $\tilde{\varphi}\in L^2(Q\times Q)$ we can write
\begin{align*}
v_n(\tilde{\varphi}) &=  \sqrt{n} \int_0^1\int_0^1 \tilde{\varphi}(x,y) \left[d\mathbb{P}_n (x,y) - dA_{1n} (x)  dA_{2n} (y)\right ] \\
&+ \sqrt{n} \varepsilon_n \int_0^1\int_0^1 \tilde{\varphi} (x,y) \left[h_n(x) - h_n(y) \right ]dQ(x)dQ(y)\\ 
&- \sqrt{n}\varepsilon_n^2 \int_0^1\int_0^1 \tilde{\varphi} (x,y) h_n(x) h_n(y) dQ(x) dQ(y), \\
\end{align*} 
and because $\tilde{\varphi} (x,y)$ is symmetric, the middle integral on the right side is zero.  The first integral, which contains the centered part of $v_n(\tilde{\varphi})$ converges in distribution to $v_{\mathbb{Q}}(\tilde{\varphi})$. Therefore, if $\varepsilon_n = O(n^{-1/4})$, we have, under $H_0$ and under $H_a$, respectively,
\begin{equation}
\label{lim_thm} 
v_n(\tilde{\varphi}) \stackrel{d}{\rightarrow}  v_{\mathbb{Q}} (\tilde{\varphi}),\quad  \text{and} \quad  v_n(\tilde{\varphi}) \stackrel{d}{\rightarrow}  v _{\mathbb{Q}} (\tilde{\varphi}) -  \langle \tilde\varphi, h\times h\rangle_{Q \times Q}.  
\end{equation}
Here, and below in this section, we use the notation 
$$\langle \varphi, \psi \rangle_{Q \times Q}=  \int_0^1\int_0^1\varphi (x,y) \psi(x,y) dQ(x)dQ(y)$$ 
for the inner product in $L^2(Q\times Q)$. We also recall that the variance of the Gaussian random variable $v_{\mathbb{Q}} (\varphi)$ is $$ Var (v_{\mathbb{Q}}(\tilde{\varphi})) = \langle \tilde \varphi,  \tilde \varphi \rangle_{Q \times Q}. $$

In the colour-blind problem one has to use functions $\mathcal{L}^*\tilde{\varphi}$. 

To make a judgement about the asymptotic power of the linear statistics, consider the distance in total variation between two Gaussian distributions with different means and equal variances. It is given by 
$$\sup_C |N_{(\mu_1,\sigma^2)}(C) - N_{(\mu_2,\sigma^2)}(C)| ,$$
where supremum  is taken over all measurable sets (or critical regions of tests) on the real line. By its definition, this measure gives the largest possible difference between the power and the level of tests among all those that can discriminate between the two distributions. Its advantage is that there is no need to specify a particular level of a test. As one can easily  see,
$$\sup_C |N_{(\mu_1,\sigma^2)}(C) - N_{(\mu_2,\sigma^2)}(C)| = 2 N_{(0,1)} \left(\frac{|\mu_1-\mu_2|}{2\sigma}\right)- 1,$$ 
where $N_{(0,1)}(x)$ denotes the standard normal distribution function. The ratio
$$T=\frac{|\mu_1-\mu_2|}{\sigma},$$
(especially when $\mu_1=0$), is often called the signal to noise ratio. The larger this ratio is, the greater the difference between the power and the  level.

With, $U_i = \max\{X_i, Y_i\}$ and $V_i = \max\{X_i, Y_i\},$ $1 \le i \le n,$ the following result holds.
\begin{prop}\label{as-opt} The statistic of asymptotically most powerful test for testing $H_0$ against the sequence of alternatives $\mathbb{A}_n$ is of the form 
	$$\mathbb{R}^s_n(h\times h) = \frac{1}{\sqrt{n}} \sum_{i=1}^n h(U_i) h(V_i) +o_P(1).$$ 
	The distance in total variation between its asymptotic distributions, under the null and under the alternatives $\mathbb{A}_n$, is
	$N_{(0,1)}(T/2), \text{where} \;T$ represents the limit of $\sqrt{n}\eps_n^2 \| h \times h \|_{Q \times Q}$.
\end{prop}
{\bf Proof.}
From \eqref{lim_thm} it follows that the asymptotic power of the test based on the linear statistic $v_n(\mathcal{L}^*\tilde{\varphi})$ is equal to $N_{(0,1)}(T_\varphi)$, where ``signal to noise ratio" is
$$T_\varphi= \frac{\langle \mathcal{L}^* \tilde\varphi, h\times h\rangle_{Q \times Q}}{ \langle \mathcal{L}^* \tilde \varphi,  \mathcal{L}^* \tilde \varphi \rangle_{Q \times Q}^{1/2}}.$$
Now note that $h\times h$ is both symmetric and passes through $\mathcal{L}^*$. Therefore, $$\langle \mathcal{L}^* \tilde\varphi, h\times h\rangle_{Q \times Q} = \langle \mathcal{L}^* \tilde\varphi, \mathcal{L}^*(h\times h)\rangle_{Q \times Q}$$ which implies that $T_\varphi$ is maximised at $\mathcal{L}^* \tilde\varphi = \mathcal{L}^*(h\times h) = h\times h$. The statistic $v_n(h\times h)$ is equal to the sum in the display formula in the proposition. On the other hand, from \eqref{fn-param-equiv} and Proposition \ref{z_n-functionals}, it follows that $\mathbb{R}^s_n(h\times h) = v_n(h\times h) + o_P(1)$. \hfill $\Box$

\begin{example}
	\label{ex4.2}
	{\rm In the case shown in Figure \ref{fig:graphs}, the alternatives are $A_{1}(x)=x, A_{2}=x^2$ and, therefore, $Q(x)=(x+x^2)/2,$ and we have  $ \displaystyle{ h(x)=  \frac{1-2x}{1+2x}}$ with no further control over $\eps$. The proposed test statistic is of the form 
\begin{equation}\label{statistic1}
 \frac{1}{\sqrt{n}} \sum_{i=1}^{n} \frac{(2 U_i -1)(2V_i -1)}{(1+2U_i)(1+2V_i)} - \sqrt{n}\left(\int_0^1 \frac{1-2x}{1+2x} dQ_n(x) \right)^2,\end{equation}
and its signal to noise ratio, for $n=400,$  is given by $\sqrt{n}\langle h, h \rangle_{Q \times Q} = 1.98$. Therefore, the power of this linear test, directed to the chosen alternatives is essentially higher than the general (not directed) Kolmogorov-- Smirnov test, cf. Figure \ref{fig:KSalt1}. However, the test may have very low or no power for many other alternatives.
}
\end{example}
%Assume that and that we would like to test $H_0$ against the sequence of alternatives $A_{1n} \times A_{2n},$  where $\displaystyle{\frac{A_{1n}}{dQ}(x) = 1 + \varepsilon_n h(x) \mbox{ and } \frac{A_{2n}}{dQ}(x)= 1 - \varepsilon_n h(x).}$ with $\varepsilon_n=n^{-1/4}$.

%Now let $Q(x)=x$ and take simpler $h(x)=2x-1.$ This implies the alternatives of the form $A_1(x) = x + 1/3(x^2-x)$ and $A_2(x)=x - 1/3(x^2-x),$  and in this case, the test statistic can be taken as  $$ \frac{1}{31.6} \sum_{i=1}^{1000} (2 U_i -1)(2V_i -1).$$ 

\vskip 0.5cm

\section{Tests based only on maxima}\label{max}

Recall that $P_n^{(2)}$ denotes the empirical distribution function of the maxima $U_i = \max\{X_i,Y_i\}$
$$P_n^{(2)}(u) = \frac{1}{n}\sum_{i=1}^n {\bf 1}_{\{U_i\le x\}},$$
which naturally should be centered  by $Q_n^2(u)$ thus leading to an empirical process
\begin{equation}\label{EP-max}
R_n^{(2)}(u) = \sqrt{n} [ P_n^{(2)}(u) - Q_n^2(u) ].
\end{equation}
Its construction, and the form of the linear functionals from $P_n^{(2)}$,
$$\int_0^1 \alpha (u) dP_n^{(2)}(u) =\frac{1}{n}\sum_{i=1}^n \alpha (U_i) $$
can prompt one to speak about ``statistics" or tests, ``based only on maxima", even though this is not an accurate expression. Indeed, the term is $Q_n^2(u)=[P^{(1)}_n (u)+P^{(2)}_n (u)]^2/4$ and so it certainly incorporates information about minima as well. Yet, the process in \eqref{EP-max} may look as the first choice to base test statistics upon and has some nice properties. For example, its covariance function has a nice structure under $H_0$
$${\rm E}_0 R_n^{(2)}(v) R_n^{(2)}(u) = Q^2(v) [1- Q(u)]^2, \; 0 \le v\le u \le 1. $$

The natural way to study the process $R_n^{(2)}$ and its function parametric version is to embed it into the colour-blind process introduced in the previous sections. First, we see that $R_n^{(2)}$ is a restriction of  $\mathbb{R}^s_n$ to rectangles $[0,u]\times[0,u]$
$$R_n^{(2)}(u)=\mathbb{R}_n(u,u)=\sqrt{n} \left[\mathbb{P}_n(u,u) - Q_n^2(u)\right].$$
This, in particular, implies that the expected value of $R_n^{(2)}(u)$ under $H_a$  is
\begin{equation}\label{expected}
{\rm E}_a R_n^{(2)}(u) = - \sqrt{n} \eps_n^2 H_n^2(u) .
\end{equation}
Now choose the functional argument of $\mathbb{R}_n(\cdot)$ as $\varphi(x,y)=\alpha(\max(x,y))$. The function $\varphi = \varphi_\alpha$ is symmetric and it belongs to $L^2(Q\times Q)$ if and only if $ \alpha\in L^2(Q^2)$ because
\begin{equation*}
%\int_0^1\int_0^1 \varphi^2(x,y) Q(dx)Q(dy) = 
\int_0^1\int_0^1  \alpha^2 (\max\{x,y\})  dQ(x)dQ(y) =    \int_0^1 \alpha^2(u) Q^2(du), \end{equation*}
or
\begin{equation*}
\langle \varphi_\alpha,\varphi_\alpha\rangle_{Q \times Q} = \langle \alpha,\alpha \rangle_{Q^2},
\end{equation*}
leading to
$$R_n^{(2)}(\alpha) = \mathbb{R}_n(\varphi_\alpha).$$
We denote by $$\mathcal{C} = \{\varphi \in L^2(Q \times Q): \varphi (x,y) = \alpha(\max\{x,y\}), \; \alpha\in L^2(Q^2), \text{and} \;\int_0^1 \alpha(u) Q^2(du)=0\}.$$
Then studying $R_n^{(2)}(\alpha),$ with $\alpha\in L^2(Q^2)$ is equivalent to studying $\mathbb{R}_n(\varphi),$ where $\varphi \in \mathcal{C}$. Hence, \eqref{process-NCB}, \eqref{pillow-def} and Proposition \ref{z_n-functionals} imply that $$R_n^{(2)}(\alpha) = v_n(\mathcal{L}^*\varphi_\alpha) + o_P(1) ,$$ and we can focus on the linear statistics appearing on the right side.  The problem of finding the optimal linear statistic requires the maximisation of a different signal to noise ratio and opens up an interesting structure.

We introduce the Radon-Nikodym derivative
$$q(x) = \frac{dH^2(x)}{dQ^2(x)}=h(x)\frac{H(x)}{Q(x)},$$
and denote by $\varphi_q(x,y)=q(\max(x,y))$. Either from \eqref{expected} or from \eqref{lim_thm} one can derive that
\begin{align*}
E_a R_n^{(2)}(\alpha) &= E_a v_n(\mathcal{L}^*\varphi_\alpha) +o(1) = \langle \mathcal{L}^* \varphi_\alpha, h \times h \rangle_{Q \times Q} + o(1) \\
&= \langle \varphi_\alpha, h\times h \rangle_{Q \times Q} +o(1).
\end{align*}
The functions of the form $(h\times h)(x,y)=h(x)h(y)$ do not belong to the class $\mathcal{C}$ (unless $h=0$), but the expression of the expected values above suggests that the projection of $h\times h$ on $\mathcal{C}$ would be useful to consider. This projection is given by the function $\varphi_q$. Indeed, 
$$\langle \varphi_\alpha, h\times h - \varphi_q \rangle_{Q \times Q} =0,$$
so that
$$\langle \varphi_\alpha, h\times h\rangle_{Q \times Q} = \langle \varphi_\alpha, \varphi_q\rangle_{Q \times Q} =\langle \alpha, q \rangle_{Q^2}.$$ 
It seems now straightforward to choose $\alpha=q$ as an optimal test statistics. However, more clarifications are needed.

The variance of $v_n(\mathcal{L}^*\varphi_\alpha)$, as we know, is equal to $\langle \mathcal{L}^*\varphi_\alpha, \mathcal{L}^*\varphi_\alpha \rangle_{Q \times Q}$. Therefore, to find the asymptotically most powerful test against the sequence of alternatives $\mathbb{A}_n$, we need to maximise the absolute value of signal to noise ratio
$$T_\alpha =  \frac{\langle \varphi_\alpha, \varphi_q\rangle_{Q \times Q} }{\langle \mathcal{L}^*\varphi_\alpha, \mathcal{L}^*\varphi_\alpha\rangle_{Q \times Q}^{1/2}} = \frac{\langle \alpha, q\rangle_{Q^2} }{\langle \mathcal{L}^*\varphi_\alpha, \mathcal{L}^*\varphi_\alpha\rangle_{Q \times Q}^{1/2}}.$$
A useful step in this direction will be to express the denominator of the above expression in terms of the inner product in $L^2(Q^2)$. It can be verified that
$$\langle \mathcal{L}^*\varphi_\alpha, \mathcal{L}^*\varphi_\alpha\rangle_{Q \times Q} = \langle \alpha,\alpha\rangle_{Q^2} - \langle \alpha, S\alpha\rangle_{Q^2},$$
where the operator $S$ is given by
$$S\alpha (x) = \alpha(x) Q(x) + 4 \int_x^1 \alpha(t) Q(dt), $$
and we need to maximise the absolute value of
$$T_\alpha = \frac{\langle \alpha, q\rangle_{Q^2} }{[ \langle \alpha,\alpha\rangle_{Q^2} - \langle \alpha, S\alpha\rangle_{Q^2}]^{1/2}}.$$

One can go into this problem, for example  -- as a problem of calculation of support function for the convex set $$\{\alpha: \langle \alpha,\alpha\rangle_{Q^2} - \langle \alpha, S\alpha \rangle_{Q^2}\leq 1\}.$$ However, it is much simpler to {\it reverse} the point of view: for a given $\alpha $, find an alternative for which $\mathbb{R}_n(\varphi_\alpha)$ will be an optimal statistic. In the previous section, this reversal will not produce a different result, whereas in the present case the maximisation in $q$ becomes simple.

\begin{prop}
Consider a cone 
$$\mathcal{M} = \{\alpha: \int_0^u \alpha(z) dQ^2 (z) \geq 0\; \text{for all}\; u>0\} .$$
If $\alpha\in\mathcal{M}$, then the power of the statistic $\mathbb{R}_n(\varphi_\alpha)$ is the largest for local alternatives \eqref{RN1} with the corresponding $q$ equal $\alpha$.

In this case, the functions $H$ and $h$ which describe the alternatives are given by
$$H(x) = \pm\sqrt {\int_0^u \alpha(z) dQ^2 (z)} \quad \text{and} \quad h(x) =  \frac{Q(x) \alpha(x)}{2 H(x)}  .$$
\end{prop}

One can call $\mathcal{M}$ the class of {\it admissible} or {\it effective} $\alpha$-s. The need to specify such a class, i.e. choose $\varphi_\alpha$ more narrowly than from the linear space $\mathcal{C}$, is visible when we try to connect $\alpha$ and $q$. This, actually gives the proof of the proposition.

{\bf Proof.} The choice of $\alpha=q$ is valid if and only if $\alpha\in \mathcal{M}$. Indeed,
$$\alpha (x) = q(x) = \frac{dH^2(x)}{dQ^2(x)}, \quad \text{or} \quad \int_0^u\alpha(x) dQ^2(x) = H^2(u)$$
and therefore the integral has to be non-negative for all $u>0$, which implies that $\alpha$ has to belong to $\mathcal{M}$. The form of $H$ and $h$ follow.
\hfill $\Box$\\  
\vskip 0.3cm

\begin{example}
{\rm For the sake of numerical comparison, we consider the same situation as in Example \ref{ex4.2}, 
	$Q=(x+x^2)/2$ while $A_1(x)=x$ and $A_2(x) =  x^2$. As we know, this leads to 
	$$h(x)=\frac{2x-1}{1+2x}, \; H(x) = -\frac{1}{2} x(1-x)\; \text {and, therefore,} \; q(x) =\frac{(x-1)(2x-1)}{(1+x)(1+2x)} .$$
	The form of the statistic resembles the form of statistics \eqref{statistic1}, but requires centering,
	$$\frac{1}{\sqrt{n}}\sum_{i=1}^n \frac{(U_i-1)(2U_i-1)}{(1+U_i)(1+2U_i)} - \sqrt{n}\int_0^1 q(z) dQ^2_n(z).$$
	The variance of this statistic %which is, certainly less then $\int_0^1 q^2(x) dQ_a^2(x)=0.0048$, 
	numerically is equal to $0.0030$, while the shift under the alternatives becomes $\sqrt{n}\,\int_0^1 q(z) dz^3 = \sqrt{n} \,0.0048$. As a result, the
	signal to noise ratio, for $n=400$, is $1.74$. Hence, the power of the linear test here is less than what it was in Example \ref{ex4.2}, although not by much.
}
\end{example}

 \vskip 0.5cm

\section{Acknowledgement}
 
The authors thank Adrian Kennedy and Joel Bancolita for their help with figures and numerical calculations. This help was wider than what eventually is included in this text.

\end{document}